\documentclass[11pt,a4paper]{article}
\usepackage[utf8]{inputenc}
\usepackage[english]{babel}
\usepackage{lmodern} 
\usepackage{amsmath}
\usepackage{amsfonts}
\usepackage{amssymb}
\usepackage{xcolor}

\newtheorem{T}{Theorem}
\newtheorem{Le}[T]{Lemma}

\newcommand{\pr}{\noindent\textbf{\textit{Proof.}} }

\author{Ir\`{e}ne Casseli
	\footnote{
			\href{mailto:irene.casseli@univ-amu.fr}{irene.casseli@univ-amu.fr}
		}
	}
\date{\small	\textit{Aix-Marseille Universit\'{e}, I2M UMR CNRS 7373, 39 Rue F. Joliot-Curie, 13453
		Marseille Cedex 13, France}} 
\title{On products of some Toeplitz operators on polyanalytic Fock spaces}

\usepackage{hyperref}
\hypersetup{
	breaklinks=true,
	colorlinks=true,                         
	linkcolor=red!60!black, 
	citecolor=cyan!90!black, 
	urlcolor=cyan!90!black,
	linkcolor=red!60!black,
	}

\begin{document}

\maketitle

	\begin{abstract}
		The purpose of this paper is to study the Sarason's problem on Fock spaces of polyanalytic functions. Namely, given two polyanalytic symbols $f$ and $g$, we establish a necessary and sufficient condition for the boundedness of some Toeplitz products $ T_{f}T_{\bar g}$ subjected to certain restriction on $f$ and $g$. We also characterize this property in terms of the Berezin transform.
	\end{abstract}~\\
	
		
		


	\section{Statement of the result}
	
	Let us begin with some historical background on the so-called  Sarason's  problem in the context of $H^2$ and $A^2$, the 	classical Hardy and Bergman spaces of the unit disk $\mathbb{D}$. Recall that, for $\varphi\in  L^2(\partial\mathbb{D})$, the Hardy space Toeplitz operator with symbol $\varphi$ is densely defined on $H^2$ by $T_\varphi(h)=P(\varphi h)$, where $P$ denotes the Riesz-Szeg\"o projection. In the same way, using again $P$ to denote the Bergman projection, the Bergman space Toeplitz operator with symbol $\varphi\in L^2(\mathbb{D})$ on $A^2$ is given by $T_\varphi(h)=P(\varphi h)$ for a suitable $h$ in $A^2$.
		
	In both situations of $H^2$ and $A^2$, it is a well known fact that a Toeplitz operator with analytic symbol $f$ is bounded if, and only if, the symbol is bounded. Moreover, in \cite{S89}, Sarason exhibited functions $f$ and $g$ in $H^2$
	such that $T_fT_{\overline{g}}$ is bounded on $H^2$ whereas at least one of these factors is unbounded ;
	this motivates the study of the boundedness of Toeplitz products involving the symbols structure. In \cite{S94}, Sarason conjectured that a necessary and sufficient condition for the product of Toeplitz $T_fT_{\overline{g}}$ to be bounded would be like the following
	\[ \sup_{z\in\mathbb{D}} \widetilde{|f|^2}(z)\widetilde{|g|^2}(z) <\infty\] where $\widetilde{u}$ is the Berezin transform of the function $u$. 
	
	Actually, the previous condition is only necessary and the conjecture fails for both the Hardy space and Bergman space of the unit disk. Counter-examples were given in \cite{N97} and \cite{APR13}. 
	However, in the context of classical Fock spaces, Cho, Park and Zhu in \cite{CPZ14} show that the Sarason's conjecture is true.  More recently Bommier-Hato, Youssfi and Zhu generalized the results obtained in \cite{CPZ14}. In \cite{BHYZ17}, they state two necessary and sufficient conditions for boundedness of the Toeplitz product $T_fT_{\overline{g}}$ in the weighted Fock space $\mathcal{F}^2_m$ of entire square-integrable functions with respect to the Gaussian measure  \[ d\lambda_m(z)= \mathrm{e}^{-|z|^{2m}},~m\geq 1.\]
	Namely, if $f$ and $g$ are non identically zero functions in $\mathcal{F}^2_m$, they show that $T_fT_{\overline{g}}$ is bounded if and only if $f=\mathrm{e}^{q}$ and $g=c\mathrm{e}^{-q}$, with $c$ a nonzero complex constant and $q$ a polynomial of degree at most $m$, if and only if the product of Berezin transforms $\widetilde{|f|^2}\widetilde{|g|^2}$ is bounded on $\mathbb{C}$.
	
	This work studies the above results in the context of Fock spaces of polyanalytic functions. We follow the approach of \cite{BHYZ17}.\\

Given $\alpha>0$, we consider the Gaussian probability measure \[ d\mu_\alpha(z)= \frac{\alpha}{\pi} \mathrm{e}^{-\alpha|z|^2} d\lambda(z) \] where $\lambda$ is the Lebesgue area measure on the complex plane. Endowed with the usual scalar product \[ \langle f,g\rangle_\alpha = \int_\mathbb{C} f\overline{g}\,d\mu_\alpha,\] the space $L^2(\mu_\alpha)=L^2(\mathbb{C},d\mu_\alpha)$ is an Hilbert space. 
For $n\in\mathbb{N}^*$, the Fock space of $n$-analytic functions $F^2_{\alpha,n}$ is the closed subspace in $L^2(\mu_\alpha)$, endowed with the norm \[  ||f||_{2,\alpha} = \bigg( \int_\mathbb{C} |f(z)|^2 d\mu_\alpha(z)\bigg)^{1/2},\] consisting of all functions $f$ satisfying $\overline{\partial}^n f=0$. 
Basic informations about polyanalytic functions can be found in the book \cite{B91}.

%
 
The reproducing kernel of the Hilbert space $F^2_{\alpha,n}$ has been computed using various method (see for instance \cite{A10}, \cite{AIM97} or  \cite{HH13}). It can be written as 
\[ K_{\alpha,n}(z,w) =  L^1_{n-1}(\alpha|z-w|^2)\mathrm{e}^{\alpha z\overline{w}}\] where $L^\beta_k$ is the generalized Laguerre polynomials\[ L^\beta_k (x)= \sum_{j=0}^{k} (-1)^j \begin{pmatrix}
k+\beta\\k-j
\end{pmatrix} \frac{x^j}{j!}.\] We also introduce the normalized kernel function 
\[k_z^{\alpha,n}=K_{\alpha,n}(.,z)/\sqrt{K_{\alpha,n}(z,z)}.\]
Moreover, the orthogonal projection $P_{\alpha,n}\colon L^2(\mathbb{C},d\mu_\alpha)\to F^2_{\alpha,n}$ is given by \[ P_{\alpha,n}f(z)= \int_\mathbb{C} K_{\alpha,n}(z,w)f(w) d\mu_\alpha(w),\] for $f\in L^2(\mathbb{C},d\mu_\alpha)$ and $ z\in\mathbb{C}$. \\\bigskip

For a linear operator $T$ on $F^2_{\alpha,n}$ define its Berezin transform (in $F^2_{\alpha,n}$)  $B_{\alpha,n}{T}$  on  $\mathbb{C}$ as \[ B_{\alpha,n}{T} (z)= \langle Tk_z^{\alpha,n},k_z^{\alpha,n}\rangle,\quad z\in \mathbb{C}. \] 
We also define the Berezin transform (in $F^2_{\alpha,n}$)  $B_{\alpha,n}\varphi$ of a function $\varphi$ which is positive and measurable on  $\mathbb{C}$ or in $L^2(\mu_\alpha)$,     by   \[ B_{\alpha,n}\varphi(z) = \langle \varphi k_z^{\alpha,n},k_z^{\alpha,n}\rangle = \int_\mathbb{C} \varphi(w) |k_z^{\alpha,n}(w)|^2 d\mu_\alpha(w),\quad z\in\mathbb{C}.\]

Moreover, given $\varphi\in L^2(\mu_{\alpha})$, the Toeplitz operator with symbol $\varphi$ is defined on a dense subset of  $F^2_{\alpha,n}$ by $T^n_\varphi(h)=P_{\alpha,n}(\varphi h)$.
\\



The aim of this paper is to prove the following result.

\begin{T} \label{T} Let $n,m,p\in\mathbb{N}^*$, $M,N\in\mathbb{N}^*$ such that $p\leq \min(m,n)$, $M\leq\min(m-p+1, n-p+1)$ and $N\leq n-p+1$. Given two functions $f\in F^2_{\alpha,M}$ and $g\in F^2_{\alpha,N}$, both non identically zero, then the following
conditions are equivalent
\renewcommand{\theenumi}{\textit{(\roman{enumi})}}
\renewcommand{\labelenumi}{\textit{(\roman{enumi})}}
\begin{enumerate} 
\item $T^m_{f}T^p_{\overline{g}}\colon F^2_{\alpha,n}\to F^2_{\alpha,n}$ is bounded ; \label{1}
\item  there exist a polynomial $q$ of degree at most  $1$ and a nonzero complex constant $c$ such that $f=\mathrm{e}^{q}$ and $g=c\mathrm{e}^{-q}$; \label{2}
\item the product  $B_{\alpha,p}({|f|^2})B_{\alpha,p}({|g|^2})$ is bounded on $\mathbb{C}$. \label{3}
\end{enumerate}
\end{T}

Note that the choice $m=n=p$ answers the question of the boundedness on $F^2_{\alpha,n}$ af a Toeplitz product $T^n_fT^n_{\overline{g}}$ with analytic symbols.\\


Henceforth, for technical convenience and without loss of generality, we deal only with the case $ \alpha = 1$. We also denote by $F^2$ the classical analytic Fock space $F^2_1$.

	\section{Preparatory results}
	
	Here, we establish preliminaries needed in the sequel. First, obviously for each $f\in F^2_{\alpha,n}$, $P_{\alpha,n}f=f$ ; we make use of this identity that played a key role in the proof of our main theorem and we call it \emph{reproduction formula}.   This formula, combined with Cauchy-Schwarz inequality, shows that the maximum order for functions in $F^2_{\alpha,n}$ is $2$. More precisely, it can be shown that \[ |f(z)|\leq   \sqrt{n} ||f||_{2,\alpha} \,\mathrm{e}^{\frac{\alpha}{2}|z|^2} \] for 
	$f\in F^2_{\alpha,n}$ and $z\in\mathbb{C}$.  \\

	Now, the following integral estimate is stated in \cite{BHYZ17} : when $m>0$, $0 \leq d \leq m$, $N > -1$, and $a \geq 0$, there is a  positive constant $C $, independent of $a$, such that \[ \int_0^{+\infty} \mathrm{e}^{-\frac{1}{2}r^{2m}+ar^d}r^N dr \leq C(1+a)^{\frac{N+1}{m}-1} \mathrm{e}^{\frac{a^2}{2}} . \]
	Here we need a special case of the latter result ($m=d=1$) in order to estimate the norm of the product operator. 		
	
	\begin{Le}  Given $a>0$ and $N\in\mathbb{N}$, define $I_N(a)$ as  \[I_N(a)= \int_0^{+\infty} r^N\mathrm{e}^{-\frac{r^2}{2}+ra}dr. \] Then, there exist a real constant $A=A(N)$ such that $I_N(a)\leq A(1+a)^N\mathrm{e}^\frac{a^2}{2}$. \label{Lestim} \end{Le}

\section{The Toeplitz product}

In this section, we first study a very special case of Toeplitz operators whose symbols take the form $\mathrm{e}^q$, where $q$ is a complex linear polynomial. This gives a sufficient condition for boundedness of the Toeplitz product. 
Subsequently, we will actually show that the condition is also necessary, by following very closely the same arguments outlined in \cite{BHYZ17}. As a result, the symbols should be an exponential of a polynomial whose degree is less than or equal to $2$.

\begin{Le} Let $f(z)=\mathrm{e}^{az}$ and $g(z)=\mathrm{e}^{-az}$ with $a\in\mathbb{C}^*$. Then,  for any $n\in\mathbb{N}^*$ and $p\leq m$, the product $T=T^m_fT^p_{\overline{g}}$ is bounded on $F^2_{n}$.\label{L1}
\end{Le}

\pr If $h$ is a polynomial in $z$ and $\overline{z}$, from Fubini’s theorem and  the reproduction formula, we obtain  
\begin{eqnarray} Th(z) & =& \int_\mathbb{C} K_m(z,v)\int_\mathbb{C}K_p(v,w)h(w)\overline{g}(w)d\mu(w)f(v)d\mu(v) \nonumber\\
& =& \int_\mathbb{C} \int_\mathbb{C}K_m(z,v) f(v)K_p(v,w)d\mu(v)\, h(w)\overline{g}(w)d\mu(w) \nonumber\\
& =& \int_\mathbb{C} f(z)K_p(z,w) h(w)\overline{g}(w)d\mu(w) \nonumber.
\end{eqnarray}
Consequently, we have
\begin{eqnarray} |Th(z)|^2\frac{\mathrm{e}^{-|z|^2}}{\pi} & \leq & \frac{1}{\pi} \bigg( \int_\mathbb{C} |K_p(z,w)|\mathrm{e}^{\text{Re}(az-\overline{aw})} |h(w)| \mathrm{e}^{-|w|^2-\frac{1}{2}|z|^2} \frac{d\lambda(w)}{\pi} \bigg)^2\nonumber\\
& =&  \bigg( \int_\mathbb{C} H_a(z,w)|h(w)| \mathrm{e}^{-\frac{1}{2}|w|^2} \frac{d\lambda(w)}{\sqrt{\pi}} \bigg)^2,\nonumber
\end{eqnarray}
where  \[ H_c(z,w)= \pi^{-1} |K_p(z,w)|\mathrm{e}^{\text{Re}\,c(z-{w})}\mathrm{e}^{-\frac{1}{2}(|z|^2+|w|^2)}\] pour $c\in\mathbb{C}$.

Now, we consider the operator $S$, formally defined on $L^2(d\lambda)$ by \[ Sh(z) = \int_\mathbb{C} H_a(z,w)h(w)d\lambda(w).\] 
We have \[  |Th(z)|^2 \frac{\mathrm{e}^{-|z|^2}}{\pi} \leq S(\pi^{-1/2}|h|\mathrm{e}^{-\frac{1}{2}|.|^2})(z)^2.\] Using the identity $||h||_2= ||\pi^{-1/2}h\mathrm{e}^{-\frac{1}{2}|.|^2}||_{L^2(d\lambda)}$ for all $h\in L^2(\mu)$, the problem of determining when $T$
would be bounded on $F^2_{n}$ reduces to the problem of determining when the operator $S$ 
is bounded on $L^2(d\lambda)$.

For each $c\in\mathbb{C}$, we set \[ H_c(z)=\int_\mathbb{C} H_c(z,w) d\lambda(w). \]
In view of Schur's test and the identity \[ \int_\mathbb{C} H_c(w,z) d\lambda(w)=H_{-c}(z), \] we conclude that the operator $S$  would be bounded on $L^2(d\lambda)$ provided that there exist a constant $C=C(c)$ such that $H_c\leq C
$ on $\mathbb{C}$. Moreover, the norm of $T$ will not exceed $C$.\\


Let $z\in\mathbb{C}$ and $c\in\mathbb{C}$. By the triangle inequality and the translation invariance of the Lebesgue measure, the following  are valid
\begin{eqnarray}
H_c(z) & = &  \frac{1}{\pi}  \int_{\mathbb{C}} \big|L^1_{p-1}(|z-w|^2)\big|  \mathrm{e}^{\text{Re }z\overline{w}}\mathrm{e}^{\text{Re}\,c(z-{w})-\frac{1}{2}(|z|^2+|w|^2)}
\lambda(w) \nonumber\\
& \leq & \frac{1}{\pi}  \int_{\mathbb{C}} \big|L^1_{p-1}(|z-w|^2)\big|  \mathrm{e}^{|c||z-{w}|-\frac{1}{2}|z-w|^2}
\lambda(w) \nonumber\\
& = & 2 \int_0^{+\infty} \big|L^1_{p-1}(r^2)\big|   \mathrm{e}^{|c|r-\frac{r^2}{2}}rdr
 \nonumber\\
& \leq &  \sum_{j=0}^{p-1}\frac{2}{j!} \begin{pmatrix}
p\\p-1-j
\end{pmatrix}  I_{2j+1}(|c|). \nonumber
\end{eqnarray}
The above inequalities, together with Lemma  \ref{Lestim}, imply that we can find positive real constants $M_1,M_2$ such that for all $c\in\mathbb{C}$,  \[ \sup_{z\in\mathbb{C}} H_c(z) \leq M_1 \mathrm{e}^{M_2|c|^2}\] with $M_2>1/2$. This gives the desired inequality, which completes the proof.


\bigskip We shall be interested here in the converse direction in the previous lemma.  We will show that the necessary condition on the polyanalytic symbols $f$ and $g$ is also a necessary condition for boundedness of the Toeplitz product if we impose some restrictions on the order of polyanalyticity of $f$ and $g$.

\begin{Le} Assume that $p\leq \min(m,n)$, $M\leq\min(m-p+1, n-p+1)$ and $N\leq n-p+1$. Given $f\in F^2_{\alpha,M}$ and $g\in F^2_{\alpha,N}$ each not identically zero such that $T=T^m_fT^p_{\overline{g}}$ is bounded on $F^2_{n}$,
then, there are a polynomial $q$ of degree at most  $1$ and a non-zero complex constant $c$ such that $f=\mathrm{e}^{q}$ and $g=c\mathrm{e}^{-q}$.
\label{L21}
\end{Le}

\pr From the Cauchy-Schwarz inequality, when $T$ is bounded, its Berezin transform $B_nT$ is bounded on the complex plane.  

Now, fix $z,a\in\mathbb{C}$ ; when $g$ is a $N-$analytic polynomial,   
\begin{eqnarray}
T^p_{\overline{g}}k^{n}_z(a) \nonumber
	 &=& \dfrac{1}{\sqrt{K_{n}(z,z)}} \overline{ \int_\mathbb{C} K_{n}(z,w) g(w) K_{p}(w,a)  d\mu(w)}\\
& =& \dfrac{1}{\sqrt{K_{n}(z,z)}} \overline{g(z) K_{p}(z,a) } \nonumber \\&=& \sqrt{\frac{p}{n}}~ \overline{g(z)} k^{p}_z(a) \nonumber
\end{eqnarray} 
where the last equality follows from the reproduction formula of $F^2_{n}$ applied to the function $gK_{p}(.,a)\in F^2_{ N+p-1}\subset F^2_{n}$ (since $N+p-1\leq n$). 
Then, the density of polyanalytic polynomials in $F^2_{N}$ ensures that the above relation   is valid also for every $g$ in $F^2_{N}$.  

Consequently when $f$ is a $M-$analytic polynomial, again applying the reproduction formula (in $F^2_{n}$ here, since $M+p-1\leq m$), we get   
  \begin{eqnarray}
Tk^{n}_z(a) 
	& = & \sqrt{\frac{p}{n}}~\overline{g(z)} \int_\mathbb{C} K_{m}(a,w)f(w) \dfrac{K_{p}(w,z)}{\sqrt{K_{p}(z,z)}}d\mu(w)\nonumber \\\nonumber
	&=& \sqrt{\frac{p}{n}}~ f(a)\overline{g(z)} k^{p}_z(a) \nonumber
\end{eqnarray}
An approximation argument then shows that the same is true given an arbitrary $f$ in $F_M^2$.

Approximating the function $f$ by polynomials and using again the reproducing formula in $F^2_{n}$, knowing that $M+p-1\leq n$ and by density, we find that
\begin{eqnarray} B_n T(z) = \frac{p}{n} f(z)\overline{g(z)}.\label{B_nT}
\end{eqnarray}  
\\

As a consequence of Liouville's theorem (see \cite{B91}, Theorem 2.5, p. 211), $fg$ must be constant as a bounded polyanalytic function. We claim that $f$ and $g$ are analytic. To see this, since neither $f$ and $g$ vanish, set $fg=c$ with  $c\in\mathbb{C}^*$. Then $f$ and $g$ are non-vanishing polyanalytic function ; thus we can write $f(z) = P(z,\overline{z})\mathrm{e}^{f_1(z)}$ and $g(z) = Q(z,\overline{z})\mathrm{e}^{g_1(z)}$
where $P$ and $Q$ are polynomials, and $f_1$ and $g_1$ are entire functions. Identifying $\mathbb{C}[z,\overline{z}]$ with $\mathbb{C}[z][\overline{z}]$,  we deduce from the identity $fg=c$ that $PQ$ must be a constant in $\mathbb{C}[z]$. 

Next, the Weierstrass factorization theorem provides that there are a complex quadratic polynomial $q(z)=a_0+a_1z+a_2z^2$  and a non-zero complex constant $c$ such that $f=\mathrm{e}^{q}$.   \\ 

We now turn to show by contradiction that $q$ is actually linear. For this purpose, assume that $a_2\neq 0$. Consider the map $S$, defined on $\mathbb{C}\times\mathbb{C}$ by  \[S(z,w) = \langle Tk^n_z,k^n_w\rangle\] which is bounded in view of the Cauchy-Schwarz inequality since $T$ is bounded.

 Again, the reproducing formula and the approximation arguments used previously yield 
\begin{eqnarray} S(z,w) 
	&=& \sqrt{\frac{p}{n}}~ \dfrac{f(w)\overline{g(z)}K_{p}(w,z)}{\sqrt{K_{p}(z,z)K_{n}(w,w)}} \nonumber\\
&= &   n^{-1}~ f(w)\overline{g(z)} 
 L^1_{p-1}( |z-w|^{2}) \mathrm{e}^{-\frac{1}{2}|z-w|^2}
\nonumber
\end{eqnarray} 
so that \begin{eqnarray}
|S(z,w)| =   \frac{|c|}{n}\,
\big|
L^1_{p-1}(|z-w|^{2}) 
 \big|\, \mathrm{e}^{-\frac{1}{2}|z-w|^2}\,\mathrm{e}^{\text{Re}(q(w)-q(z))}
\nonumber.
\end{eqnarray} For sufficiently large $t>0$, we have $L^1_{p-1}(t^2|a_2|^2)\neq 0$. Taking $z=r\in\mathbb{R}_+$ and $w=r+t\overline{a_2}$, it follows that there exists a real constant $A=A(n,c,a_1,a_2)$ such that  \begin{eqnarray}
|S(r,r+\overline{a_2})| = A \mathrm{e}^{2|a_2|^2r}.
\nonumber
\end{eqnarray} We reach a contradiction with the boundedness of $S$ when $a_2\neq 0$.


\bigskip To sum up, we have proved the following statement which corresponds to the equivalence between \ref{1} and \ref{2} of our main theorem. 

\begin{T} Let $n,m,p\in\mathbb{N}^*$, $M,N\in\mathbb{N}^*$ such that $p\leq \min(m,n)$, $M\leq\min(m-p+1, n-p+1)$ and $N\leq n-p+1$. If $f\in F^2_{\alpha,M}$ and $g\in F^2_{\alpha,N}$ each non identically zero, then the Toeplitz product $T^m_{f}T^p_{\overline{g}}$ 
is bounded on $F^2_{n}$ if and only if  $f=\mathrm{e}^{q}$ and  $g=c\mathrm{e}^{-q}$ where $q$ is a complex linear polynomial and  $c$ is a non-zero complex constant. \label{Tintermed}
\end{T}

\section{Sarason's conjecture}

In what follows, we provide a solution to Sarason's problem for some Toeplitz products with polyanalytic symbols in the Fock space of polyanalytic functions. Namely, thanks to Theorem \ref{Tintermed} of the above section, it becomes clear that Sarason's conjecture turns out to be true for polyanalytic Fock spaces setting.  We will show that condition \ref{3} of Theorem \ref{T} stated in introduction is equivalent to conditions \ref{1} and \ref{2} by separating it into two lemmas. Again,  our proof follows the same arguments stated in \cite{BHYZ17}.

 We first show that Berezin transforms of the square of the modulus of any polyanalytic function $h$  pointwise controls $|h|^2$.

 \begin{Le} \label{Majoration} Suppose that $m,n\in\mathbb{N}^*$ and $h\in F^2_n$. 
 Then  \[|h|^2\leq \frac{{m+n-1}}{m}  B_{m}({|h|^2})\] on $\mathbb{C}$.
\end{Le}

\pr
If  $h$ is a polyanalytic polynomial in $F^2_n$, then by virtue of the reproduction formula at the point $z\in\mathbb{C}$, it follows that
 \[ h(z) K_m(z,z) = \int_\mathbb{C} K_{m+n-1}(z,w)h(w)K_m(w,z) d\mu(w). \]
 
This equality, combined with an approximation argument and   the Cauchy-Schwarz inequality, implies that
 \begin{eqnarray}
 |h(z)|^2 &\leq&   \Bigg(\int_\mathbb{C} |K_{m+n-1}(z,w)h(w)\frac{K_m(w,z)}{K_m(z,z)}| d\mu(w) \Bigg)^2\nonumber\\
 	&\leq&   \int_\mathbb{C} \bigg|\frac{K_{m+n-1}(z,w)}{\sqrt{K_m(z,z)}}\bigg|^2 d\mu(w) \int_\mathbb{C} |h(w)|^2|k^m_z(w)|^2 d\mu(w)\nonumber\\
 	 & = &   \frac{{m+n-1}}{m}   B_{m}({|h|^2}) \nonumber .
 \end{eqnarray}


We keep throughout the rest of the paper the hypotheses of our main theorem, that is $f\in F^2_{\alpha,M}$ and $g\in F^2_{\alpha,N}$ are non-identically 
zero where $n,m,p,M,N\in\mathbb{N}^*$ such that $p\leq \min(m,n)$, $M\leq\min(m-p+1, n-p+1)$ and 
$N\leq n-p+1$. 
 As a consequence of the previous lemma, the following result can be established.
 
\begin{Le} \label{Le} If  $B_{\alpha,p}({|f|^2})B_{\alpha,p}({|g|^2})$ is bounded on $\mathbb{C}$, then the Toeplitz product $T=T^m_{f}T^p_{\overline{g}}$ is bounded on $F^2_{n}$.
\end{Le}

\pr Applying of Lemma \ref{Majoration} shows that, when $B_{\alpha,p}({|f|^2})B_{\alpha,p}({|g|^2})$ is bounded on $\mathbb{C}$, the same is true for $fg$ ; the arguments given in the proof of Lemma \ref{L21} ensure that there exists a non-zero complex constant $c$ and a complex quadratic polynomial $q(z)=a_0+a_1z+a_2z^2$ and a non-zero complex constant $c$ such that $f=\mathrm{e}^{q}$ and $g=c\mathrm{e}^{-q}$.

As in for the previous proof, let us assume that $a_2\neq 0$ and show that this leads to a contradiction. 
  Define a map $B$ on $\mathbb{C}$ by setting 
\[B = |f|^2B_{p}({|g|^2}).\]
 This map is bounded in view of Lemma \ref{Majoration}. Now, for every $x\in\mathbb{R}_+$, 
 \begin{eqnarray}
  |B(x)|^2
   & = &  \mathrm{e}^{2\text{Re}q(x)}  \int_\mathbb{C}    
   |c|^2\mathrm{e}^{-2\text{Re}q(w)} 
   \frac{\big| L^1_{p-1}(|x-w|^{2}) 
   	\big|^2}{p}
	\, \mathrm{e}^{2\text{Re}x\overline{w}-|x|^2} \frac{\mathrm{e}^{-|w|^2} }{\pi}
	d\lambda(w) \nonumber\\
	& = &  \frac{|c|^2}{p\pi} \int_\mathbb{C}
	\mathrm{e}^{2\text{Re}(q(x)-q(w))} 	
   \big| L^1_{p-1}(|x-w|^{2}) 
   \big|^2
   \mathrm{e}^{-|x-w|^2}
	d\lambda(w)\nonumber\\
	&\geq &  \frac{|c|^2}{p\pi} \int_\mathbb{C}
	\mathrm{e}^{2\text{Re}(a_2(x^2-w^2))} 	
   \big| L^1_{p-1}(|x-w|^{2}) 
   \big|^2
   \mathrm{e}^{-|x-w|^2-2|a_1(x-w)|}
	d\lambda(w).\nonumber
\end{eqnarray}


Since $ L^1_{p-1}$ is a polynomial, one can find strictly positive real constants $R$ and $M$ such
 that \[| L^1_{p-1}(|\zeta|^{2}) \big|^2\geq M \]   for all $\zeta\in\mathbb{C}$ with $|\zeta|\geq R $. 

Set $a_2=|a_2|\mathrm{e}^{\mathrm{i}\beta}$. Inserting the previous estimate for the integrand in the last displayed inequalities and using a suitable 
change of variables, we obtain 
\begin{eqnarray}
  |B(x)|^2 
	&\geq &  \frac{|c|^2}{p\pi} \int_\mathbb{C}
	\mathrm{e}^{2\text{Re}(a_2(\zeta^2+2\zeta x))} 	
   \big| L^1_{p-1}(|\zeta|^{2}) 
   \big|^2
   \mathrm{e}^{-|\zeta|^2-2|a_1\zeta|}
	d\lambda(\zeta)\nonumber\\
	&\geq & \frac{M|c|^2}{p\pi}
	 \int_R^{+\infty} \int_{|\theta+\beta|<\frac{\pi}{4}} 	 	
   \mathrm{e}^{-(1+2|a_2|)r^2-2|a_1|r}
	rd\theta dr \,\mathrm{e}^{2\sqrt{2}R|a_2|x} .\nonumber
\end{eqnarray} 

Consequently, there exist real constants 
$A_1=A_1(n,c,a_1,a_2)$ and $A_2=A_2(n,a_2)$, with $A_2>0$, such that for all $x>0$, we have \begin{eqnarray}
|B(x)|^2\geq  A_1 \mathrm{e}^{A_2x}.
\nonumber
\end{eqnarray} 
This yields a contradiction since $B$ should be bounded.\\

Finally, we turn to the converse of the latter lemma, that is :

\begin{Le} \label{Le13} 
Assume that $T=T^m_{f}T^p_{\overline{g}}$ is bounded on $F^2_{n}$ ;
 then $B_{\alpha,p}({|f|^2})B_{\alpha,p}({|g|^2})$
 is a bounded map on  $\mathbb{C}$.
\end{Le}

\pr If $T$ is bounded, given the equalities already proven in the proof of Lemma \ref{L21}, we claim that for each  
$z\in\mathbb{C}$,
\begin{eqnarray}
 \langle Tk^n_z, T k^n_z\rangle &=& \frac{p}{n} |g(z)|^2 \langle fk^p_z, fk^p_z\rangle \nonumber\\
&=& \frac{p}{n} |g(z)|^2 \int_\mathbb{C} |f(w)|^2 |k^p_z(w)|^2 d\mu(w)\nonumber\\
&=& \frac{p}{n} |g(z)|^2 B_{p}({|f(z)|^2}) .\nonumber
\end{eqnarray}
By the Cauchy-Schwarz, $|g|^2 B_{p}({|f|^2})$  must be bounded. 

Moreover, we have $(T^m_fT^p_{\overline{g}})^*=T^m_gT^p_{\overline{f}}$. It is a consequence of the Fubini theorem together with an approximation argument.
By symmetry, since $T^*$ is bounded, we get also 
that  $|f|^2 \widetilde{|g|^2}$ is a bounded map. 

But, again using the proof of Lemma \ref{L21}, once $T$ is continuous, the 
product $fg$ is constant and the desired result follows namely 
$B_{\alpha,p}({|f|^2})B_{\alpha,p}({|g|^2})$ is bounded.\\

\section{Proof of the main result}

Finally, we prove Theorem \ref{T}. That \ref{1} and \ref{2} are equivalent follows from Theorem \ref{Tintermed}. To prove that \ref{1} implies \ref{3} we use Lemma \ref{Le13}. To show that \ref{3} implies \ref{1} we apply Lemma, and hence the proof is complete.

\section*{Concluding remarks}
\begin{enumerate}
\item It would be of interest to prove Theorem \ref{T} without restriction on degree of polyanalyticity of $f$ and $g$.
\item It would be also interesting to carry out this study for generalized Fock spaces of polyanalytics.
\end{enumerate}

\section*{Acknowledgements}

This work is part of my PhD thesis. I am grateful to my advisers Professors S. Rigat and E. H. Youssfi for suggesting
this problem and helpful discussions.




\begin{thebibliography}{10}
	
	\bibitem{A10} {L. D. Abreu, 
		\textit{ On the structure of Gabor and super Gabor spaces}, Monatsh. Math., 161, 237-253 (2010).} 		
	
	\bibitem{APR13} 
	A. Aleman, S. Pott, M. C. Reguera,   \textit{Sarason Conjecture on the Bergman space}, preprint (2013) available at \href{https://arxiv.org/abs/1304.1750}{https://arxiv.org/abs/1304.1750}. 
	
	\bibitem{AIM97}
	{N. Askour, A. Intissar, Z. Mouayn, \textit{{Espaces de Bargmann généralisés et formules explicites pour leurs noyaux reproduisants}}, C. R. Acad. Sci. Paris Sér. I Math., 707–712 325 (1997).}	
	
	\bibitem{B91}
	{
		M. B. Balk, \textit{Polyanalytic Functions}, Akad. Verlag, Berlin (1991).}	
	
	
	
	\bibitem{BHYZ17}
		{H. Bommier-Hato, E. H. Youssfi, K. Zhu, \textit{Sarason’s Toeplitz product problem for a class of
				Fock spaces}, Bull. Sci. math., 141,  408–442 (2017).}
		
		\bibitem{CPZ14}		
			{	H.R. Cho, J.D. Park, K. Zhu, \textit{Products of Toeplitz operators on the Fock space}, Proc. Am. Math.
				Soc., 142,  2483–2489 (2014).}
			
			
			\bibitem{HH13}
			{	A. Haimi, H. Hendenmalm,
				\textit{The polyanalytic Ginibre ensembles}, J. Stat. Phys., 153, 1,  10–47 (2013). }
			
			\bibitem{N97}
			F. Nazarov,
			\textit{A counterexample to Sarason's conjecture}, preprint (1997), available at \href{ http://users.math.msu.edu/users/fedja/prepr.html}{ http://users.math.msu.edu/users/fedja/prepr.html}.
			
			\bibitem{S89}
			{D. Sarason, 
				\textit{Exposed points in $H^1$} I, The Gohberg anniversary collection, Vol. II,
				pp. 485–496, Oper. Theory Adv. Appl., 41, Birkhäuser, Basel, (1989).  }
			
			\bibitem{S94}  
			{D. Sarason,
				\textit{Products of Toeplitz operators}, in \textit{Linear and Complex Analysis Problem Book 3}, Part I (V. P. Havin and N. K. Nikolskii, Eds.), part of Lecture Notes Math., Vol. 1573, pp. 318-319, Springer-Verlag, Berlin-New York, (1994).}
	

%
%
%
%
%
%
%
%
%
%
%
%

		
		

\end{thebibliography}
\end{document}